\documentclass[twoside,11pt]{article}

\usepackage{graphicx, subfigure,pstricks, pst-node, color}
\usepackage[top=1in,bottom=1in,left=1in,right=1in]{geometry}
\usepackage[sort&compress]{natbib} 
\usepackage{amsfonts, amsmath, amssymb, amsthm}
\usepackage{psfrag}
\usepackage{subfigure}
\usepackage{flafter}
\usepackage[section]{placeins}
\usepackage{color}
\definecolor{darkred}{RGB}{100,0,0}
\definecolor{darkgreen}{RGB}{0,100,0}
\definecolor{darkblue}{RGB}{0,0,150}

\usepackage{hyperref}
\hypersetup{colorlinks=true, linkcolor=darkred, citecolor=darkgreen, urlcolor=darkblue}
\usepackage{url}

\newtheorem{thm}{Theorem}
\newtheorem{prp}{Proposition}
\newtheorem{lem}{Lemma}

\def\beq{\begin{equation}}
\def\eeq{\end{equation}}
\def\beqn{\begin{eqnarray*}}
\def\eeqn{\end{eqnarray*}}
\def\bitem{\begin{itemize}}
\def\eitem{\end{itemize}}
\def\benum{\begin{enumerate}}
\def\eenum{\end{enumerate}}
\def\bmult{\begin{multline*}}
\def\emult{\end{multline*}}
\def\bcenter{\begin{center}}
\def\ecenter{\end{center}}

\newcommand{\thmref}[1]{Theorem~\ref{thm:#1}}
\newcommand{\prpref}[1]{Proposition~\ref{prp:#1}}

\newcommand{\lemref}[1]{Lemma~\ref{lem:#1}}
\newcommand{\secref}[1]{Section~\ref{sec:#1}}

\DeclareMathOperator*{\argmax}{arg\, max}
\DeclareMathOperator*{\argmin}{arg\, min}
\def\cJ{\mathcal{J}}

\def\cN{\mathcal{N}}

\def\bA{\boldsymbol{A}}

\def\bI{\boldsymbol{I}}

\def\bP{\boldsymbol{P}}
\def\bQ{\boldsymbol{Q}}

\def\bX{\boldsymbol{X}}

\def\bu{\boldsymbol{u}}

\def\bx{\boldsymbol{x}}
\def\by{\boldsymbol{y}}
\def\bz{\boldsymbol{z}}

\newcommand{\bbeta}{{\boldsymbol\beta}}

\newcommand{\bxi}{{\boldsymbol\xi}}

\newcommand\bSigma{{\boldsymbol\Sigma}}

\def\bbR{\mathbb{R}}

\newcommand{\E}{\operatorname{\mathbb{E}}}
\renewcommand{\P}{\operatorname{\mathbb{P}}}
\newcommand{\Var}{\operatorname{Var}}

\newcommand{\pr}[1]{\mathbb{P}\left(#1\right)}

\newcommand{\1}{{\rm 1}\kern-0.24em{\rm I}}

\newcommand{\<}{\langle}
\renewcommand{\>}{\rangle}

\def\eps{\varepsilon}

\def\bbhat{\widehat{\bbeta}}

\def\smax{\overline{s}}
\def\lmin{\nu}
\def\Id{\mathbf{I}}
\def\vsign{\mathrm{sign}}

\def\btrue{\bbeta_\star}
\def\strue{s_\star}
\def\Jtrue{J_\star}

\def\Jhat{\widehat{J}}
\def\es{\eps}
\def\Jmap{\Jhat_{\rm map}}
\def\bbmap{\bbhat_{\rm map}}
\def\bbpm{\bbhat_{\rm mean}}

\begin{document}

\title{Variable Selection with Exponential Weights \\
and $\ell_0$-Penalization}
\author{
Ery Arias-Castro%
%\footnote{Department of Mathematics, University of California, San Diego
%\{\href{mailto:eariasca@ucsd.edu}{eariasca@ucsd.edu}\}}
\ and
Karim Lounici%
%\footnote{School of Mathematics, Georgia Institute of Technology
%\{\href{mailto:klounici@math.gatech.edu}{klounici@math.gatech.edu}\}}
}
\date{
University of California, San Diego
\ and \
Georgia Institute of Technology
}
\maketitle

\begin{abstract}
In the context of a linear model with a sparse coefficient vector, exponential weights methods have been shown to be achieve oracle inequalities for prediction.  We show that such methods also succeed at variable selection and estimation under the necessary identifiability condition on the design matrix, instead of much stronger assumptions required by other methods such as the Lasso or the Dantzig Selector.  The same analysis yields consistency results for Bayesian methods and BIC-type variable selection under similar conditions.
\medskip

%\noindent {\bf MSC 2010:}

\noindent {\bf Keywords:} Variable selection, model selection, sparse linear model, exponential weights, Gibbs sampler, identifiability condition.

\end{abstract}

\section{Introduction}
\label{sec:intro}

Consider the standard linear regression model:
\beq \label{model}
\by = \bX \btrue + \bz,
\eeq
where $\by \in \bbR^n$ is the response vector; $\bX \in \bbR^{n \times p}$ is the regression (or design) matrix, assumed to have normalized columns; $\btrue \in \bbR^p$ is the coefficient vector; and $\bz \in \bbR^n$ is white Gaussian noise, i.e., $\bz \sim \cN(0, \sigma^2\bI_n)$. {As in general the model \eqref{model} is not identifiable, we let $\btrue$ denote one of the coefficient vectors such that $\bX \bbeta = \E(\by)$ of minimal support size.  Then $\Jtrue$ and $\strue$ denote the support and support size of $\btrue$.
We are most interested in the case where the coefficient vector is sparse, meaning $\strue$ is much smaller than $p$.} As usual, we want to perform inference based on the design matrix $\bX$ and the response vector $\by$.  The three main inference problems are:
\bitem \setlength{\itemsep}{0in}
\item {\em Prediction:} estimate the mean response vector $\bX \btrue$;
\item {\em Estimation:} estimate the coefficient vector $\btrue$;
\item {\em Support recovery:} estimate the support $\Jtrue$.
% = \{j\,:\, \bbeta_{\star,j}\neq 0\}$ of the coefficient vector $\btrue$.
\eitem

These problems are not always differentiated and often referred to
jointly as {\em variable/model selection} in the statistics
literature, and {\em feature selection} in the machine learning
literature.  Being central to statistics, a large number of papers
address these problems.  We review the literature with particular
emphasis on papers that advanced the theory of model selection.
For penalized regression, we find \citep{MR1466682}, who provides
necessary conditions and sufficient conditions under which the
AIC/Mallows' $C_p$ criteria and the BIC criteria are consistent. For
example, AIC/Mallows' $C_p$ are consistent when there is a unique
$\bbeta$ such that $\E(\by) = \bX \bbeta$, and this $\bbeta$
has a support of fixed size as $n, p \to \infty$.  Also, BIC is
consistent when the dimension $p$ is fixed and the model is
identifiable --- a condition that appears to be missing in that
paper. BIC was recently shown in \citep{MR2443189} to be consistent
when the model is identifiable, $p = O(n^a)$ with $a < 1/2$ and
the true coefficient vector has a support of fixed size as $n, p \to
\infty$.  They also propose an extended BIC for when $a$ is larger.
Assuming the size of the support of $\btrue$ is known, \cite{raskutti2009minimax} establish {\em prediction} and {\em estimation} performance bounds for best subset selection, and obtain information bounds for these problems.
Relaxing to the $\ell_1$-norm penalty, the Lasso and the closely related
Dantzig Selector were shown to be consistent when the design matrix
satisfies a restricted isometric property (RIP) or has column
vectors with low coherence; see
\citep{BRT07,MY06,B07,L08,BTW07,ZY06,MR2543688,MR2382644} among
others.  With a carefully chosen nonconcave penalty,
\citep{MR2065194} shows that consistent variable selection is
possible when $p = O(n^{1/3})$.  This condition on $p$ was weaken in
the follow-up paper \citep{MR2849368}, though with an additional
restriction on the coherence (Condition (16) there).  The strongest
results in that line of work seem to appear in \citep{MR2604701},
which suggests a minimax concave penalty that leads to consistent
variable selection under much weaker assumptions.
The classical forward stepwise selection, also known as orthogonal
matching pursuit, which is shown in \citep{CaiOMP} to enable
variable selection under an assumption of low coherence on the
design matrix. Screening was studied in \citep{fan2008sure} in the
ultrahigh dimensional setting, assuming the design is random.  A
combination of screening and penalized regression is explored in
\citep{ji2010ups,jin2012optimality}, with asymptotic optimality when
the Gram matrix $\bX^\top\bX/n$ is (mildly) sparse.

A distinct line of research is the implementation of
$\ell_0$-penalized regression via exponential weights
\citep{Ca2004,Yan04,DalTsy2007,Gir07,JuRiTsy2006,L07,SalDal11}. This
methodology, which has precedents in the Bayesian literature on
model selection \citep{MR2000752}, has the potential of striking a
good  compromise between statistical accuracy and computational
complexity.  While computational tractability has only been
demonstrated in simulations, a number of sharp statistical results
exist for the {\em prediction} problem.  In particular,
\citep{AL11,MR2816337} propose exponential weights procedures that
achieve sharp sparsity oracle inequalities with no assumptions of
the design matrix $\bX$. Note that there exists no result in the literature concerning the problems of estimation and support recovery with an exponential weights approach.  For a recent survey of the exponential weights literature, see \citep{RigTsy11}.

Our contribution is the following.  We establish performance bounds
for the version of exponential weights studied in \citep{AL11} for
the three main inference problems of {\em prediction}, {\em estimation} and {\em support recovery}. The methodoly developed in the present paper is new and brings novel and interesting results to the sparse regression literature. The
main feature of this methodology is that it only requires
comparatively almost minimum assumptions on the design matrix $\bX$.  In particular, for {\em estimation} and {\em support recovery}, the conditions are slightly stronger than identifiability.  Moreover, when the size of support is known, the exponential weights method is consistent under the minimum identifiability condition as long as the nonzero coefficients are large enough, close in magnitude to what is required by any method, in particular matching the performance of best subset selection \citep{raskutti2009minimax}.  See also \citep{verzelen,candes-davenport,zhang}.  An important by-product of our analysis are
consistency results for BIC-type methods, i.e., variable selection
with $\ell_0$-penalty, under similar conditions, extending the results of \cite{MR2443189}.

The rest of the paper is organized as follows.  In \secref{main}, we
describe in detail the methodology and state the main results. We
also state similar results for variable selection with
$\ell_0$-penalty and for the Bayesian model selection method of
\citep{MR2000752}.  In \secref{comparison}, we compare the results
we obtained for exponential weights with those established for
other methods, in particular the Lasso and {\sc mc+}.  In \secref{num}, we briefly
discuss the algorithmic implementation of exponential weights,
and show the result of some simple numerical experiments comparing
exponential weights with other popular variable selection techniques in the literature. In \secref{discussion}, we
discuss our results in the light of recent information bounds for model selection.
%obtained for the problems of {\em estimation} and {\em support
%recovery}.
The proofs of our main results are in \secref{proofs}.

\section{Main results}
\label{sec:main}

We consider the version of exponential weights studied in
\citep{AL11}, shown there to enjoy optimal oracle performance for
the prediction problem. {The procedure puts a sparsity prior on the coefficient vector and selects the estimates using the posterior distribution.
We obtain a new {\em prediction} performance bound which is based on balancing the sparsity level and the size of the least squares residuals.  The result does not assume any conditions on the design matrix.  The task of {\em support recovery}, to be amenable, necessitates additional assumptions.  We show that under near-identifiablity conditions on the design matrix, the posterior concentrates on the correct subset of nonzero components with overwhelming probability, provided that these coefficients are sufficiently large --- somewhat larger than the noise level.  This immediately implies that the maximum a posteriori (MAP) is consistent.  We then derive {\em estimation} performance guarantees in Euclidean norm and $l_{\infty}$-norm for the maximum a posteriori and posterior mean.}

Throughout, we assume the noise variance $\sigma^2$ is known.  We also assume that $p \ge n$ and remark that similar results hold when $n \ge p$, with $p$ replaced by $n$ in the bounds.
%Implementation issues are postponed to \secref{num}.

We use some standard notation.  For any $\bu=(\bu_1,\cdots,\bu_d)^\top\in \bbR^d$ with $d\ge 1$ and $q \ge 1$,  we define
$$
\|\bu\|_q = \left(\sum_{j=1}^d \bu_j\right)^{1/q},\quad \|\bu\|_\infty = \max_{1\le j \le d} |\bu_j|.
$$
Without loss of generality, we assume from now on that the predictors are normalized in the sense that \beq \label{normalized}
\frac{1}{\sqrt{n}}\|\bX_j\|_2 = 1, \text{ for all $1\le j \le p$}.
\eeq
For a subset $J \subset [p] := \{1, \dots, p\}$, let $\bX_J = [\bX_j, j \in J] \in \bbR^{n \times |J|}$, where $\bX_j$ denotes the $j$th column vector of $\bX$.
For a subset $J \subset [p]$, let $M_J$ be the linear span of $\{\bX_j, j\in J\}$ and let $\bP_J$ be the orthogonal projection onto $M_J$. Then, $\bP_{J}^\perp := \bI_n - \bP_J$ is the orthogonal projection onto $M_J^\perp$.  We say that a vector is $s$-sparse if its support is of size $s$.

\subsection{Exponential weights}
\label{sec:expo}

We start with the definition of a sparsity prior on the subsets of $[p]$, which favors subsets with small support.  This leads to a pseudo-posterior, which is used in turn to define various exponential weights estimators.

\bitem
\item \emph{The prior $\pi$.} Fix an upper bound $\smax \geq 1$ on the support size, and a sparsity parameter $\lambda>0$.  The prior chooses the subset $J \subset [p]$ with probability
\begin{align}\label{prior-sparse}
\pi(J) \propto \binom{p}{|J|}^{-1} e^{-\lambda|J|} \, \1_{\{|J| \le
\smax\}}.
\end{align}
%where $C(\smax)$ is the renormalization constant.
\item \emph{The posterior $\Pi$.} Given that the noise is assumed i.i.d.~Gaussian with variance $\sigma^2$, given a subset of variables $J \subset [p]$, the coefficient vector that maximizes the likelihood is the least squares estimate $\bbhat_J$ with a maximum proportional to $\exp\left(-\|\bP_J^\perp(\by)\|_2^2/(2 \sigma^2)\right)$.  In light of this, we define the following pseudo-posterior, which chooses $J \subset [p]$ with probability
\begin{align}\label{posterior-sparse}
\Pi(J) \propto \pi(J) \, \exp\left(-\frac{\|\bP_J^\perp(\by)\|_2^2}{2
\sigma^2}\right).  %, \, \forall J \subset [p].
\end{align}
\eitem
The prior $\pi$ enforces sparsity and focuses on subsets of size not exceeding $\smax$. Without additional knowledge, we shall take $\smax = p$.
%[n/2] \wedge p$ since, in view of \lemref{id}, the model \eqref{model} is not identifiable for coefficient vectors of support size exceeding $[n/2]$.
%
The exponential factor in $\|\bP_{J}^\perp(\by)\|_2^2$ in the posterior enforces fidelity to the observations.  Note that $\Pi$ is not a true posterior because no prior is assumed for $\btrue$; we elaborate on this point in \secref{bayes}.  The variance term $2 \sigma^2$ corresponds to the temperature $T$ in a standard Gibbs distribution.  We will calibrate the procedure via the sparsity exponent $\lambda$ in \eqref{prior-sparse}, though we could have done so via the temperature as well.  Remember that we assume that $\sigma^2$ is known.  When the variance is unknown, we can replace it with a consistent estimator $\hat\sigma^2$.

Based on the pseudo-prior $\Pi$, it is natural to consider the maximum a posteriori (MAP) support estimate, defined as
\beq \label{map}
\Jmap = \argmax_J \, \Pi(J).
\eeq
This leads to considering the MAP coefficient estimate.  For any $J \subset [p]$, let $\bbhat_J$ denote the the least squares coefficient vector for the sub-model $(\bX_J,\by)$ with minimum Euclidean norm --- so that $\bbhat_J$ is unique even when the columns of $\bX_J$ are linearly dependent.  When the  columns of $\bX_J$ are linearly independent, the standard formula applies
\beq \label{LS}
\bbhat_J =  (\bX_J^\top \bX_J)^{-1}\bX_J^{\top}\by.
\eeq
Note that $\Id(s)$ guarantees that \eqref{LS} holds when $|J| \leq s$.  The MAP coefficient estimate is then defined as $\bbmap = \bbhat_{\Jmap}$.

We found that the MAP is not as stable as the posterior mean
\beq \label{PM}
\bbpm = \sum_J \Pi(J) \bbhat_J.
\eeq
We establish results for both of them.
%Though in our simulations the posterior mean performs better than the MAP, the theoretical performance bounds we obtain are better for the MAP than the posterior mean.

\subsection{Prediction}

We establish a new sparsity oracle inequality for the prediction problem.  We show that, in terms of prediction performance, the maximum a posteriori and posterior mean come within a log factor of that of the oracle estimator $\bbhat_{\Jtrue}$:
\[
\|\bX \bbhat_{\Jtrue} -\bX \btrue\|_2 = \|\bP_{\Jtrue} \bz \|_2 = O_P(\sigma \sqrt{\strue}).
\]

\begin{thm} \label{thm:prediction}
Consider a design matrix $\bX$ with $p \ge n$ and normalized column vectors \eqref{normalized}.
Assume $\lambda = (62 + 12c) \log p$ for some $c > 0$.  Then with probability at least $1 - p^{-c}$,
\beq \label{pred}
\|\bX \bbmap -\bX \btrue\|_2 \le \sigma \sqrt{8 \strue \lambda}
\quad \text{and} \quad
\|\bX \bbpm -\bX \btrue\|_2 \le \sigma \sqrt{12 \strue\lambda}.
\eeq
\end{thm}
Note that here, and anywhere else in the paper, what is true of $\bbmap$ is true of $\bbhat_J$ for any $J$ such that $\Pi(J) \ge \Pi(\Jtrue)$.

In \citep{AL11}, a similar sparsity oracle inequality is established in expectation using the approach by Stein's Lemma from \citep{LB06}.  Here, we use instead the concentration property of the posterior $\Pi$ and show that the oracle inequality also holds true in probability. Note that \cite{AL11} also established an oracle inequality in probability for a different exponential weights procedure that requires the knowledge of $\|\btrue\|_1$.  Our result constitutes an improvement since we do not require such knowledge.

\subsection{Concentration of the posterior and support recovery}

Our performance bounds for support recovery rely, as they should, on concentration properties of the posterior $\Pi$.  We first prove that, without any condition on the design matrix $\bX$, the posterior $\Pi$ concentrates on subsets of small size.

\begin{prp}\label{prp:support}
Consider a design matrix $\bX$ with $p \ge n$ and normalized column vectors \eqref{normalized}.
For some $\es > 0$ and $c \ge 1$, take
\beq \label{lambda}
\lambda = \frac{1+\es}{\es} (23+5c) \log p.
\eeq
Then, with probability at least $1-2p^{-c}$,
$\Pi(J) < \Pi(\Jtrue)$ for all $J \subset [p]$ such that $|J| > (1+\es) \strue$, and in fact
\beq \label{support}
\Pi\left(J: |J| > (1+\es) \strue \right) \le 4 p^{-c} \, \Pi(\Jtrue).
\eeq
%\item For some $\es, c, M > 0$, take
%\beq \label{lambda-2}
%\lambda = \frac{(1+\es) (23+5c) + M}{\es} \log p.
%\eeq
%Then, with probability at least $1-2p^{-c}$,
%\begin{align*}
%\sum_{J\,:\,  |J| > (1+\es) \strue } |J|^m \frac{\Pi\left(J\right)}{\Pi(\Jtrue)} \le 3p^{-c},\quad \forall 0\leq m \leq M.
%\end{align*}
%\eenum
\end{prp}

\subsubsection{Identifiability}

Actual support recovery requires some additional conditions, the bare minimum being that the model is identifiable.
\begin{description}
\item[Condition $\Id(s)$:]
For any subset $J \subset \{1, \dots, p\}$ of size $|J| \leq s$, the submatrix $\bX_J$ is full-rank.
\end{description}
This condition characterizes the identifiability of the model as stated in the following simple result.

\begin{lem} \label{lem:id}
Assuming $\btrue\in \bbR^p$ is $\strue$-sparse, it is identifiable if, and only if, $\Id(2\strue)$ is satisfied.
\end{lem}
%The proof is straightforward and left to the reader.

In this paper, we establish that exponential weights, and also
$\ell_0$-penalized variable selection, allow for support recovery and estimation under the condition
$\Id((2+\es)\strue)$ for any $\es > 0$ fixed, as long as the non-zero entries of the coefficient vector are sufficiently large.  In fact, $\Id(2\strue)$ suffices when $\strue$ is known.
%(Note that $\smax$ is given).

%For any $\bu\in \bbR^p$ and $J\in [p]$, we denote by $\bu_J$ the vector of $\bbR^p$ obtained by %keeping the components $\bu_j$ of $\bu$ with index $j
%\in J$ and putting the remaining ones equal to $0$.

While $\Id(s)$ is qualitative, results on estimation and support recovery necessarily require a quantitative measure of correlation in the covariates. The following quantity appears in the performance bounds we derive for exponential weights and related methods:  for any integer $s \geq 1$, define
\beq \label{lmin}
%\lmin_{s} = \min\left\{\frac{1}{\sqrt{n}} \|\bX_J \, \bu\|_2 : J \subset [p] \text{ s.t. } |J|\leq s; \bu\in \bbR^{|J|} \text{ s.t. } \|\bu\|_2 = 1\right\}.
\lmin_{s} = \min_{J \subset [p]\,:\, |J|\leq s} \ \min_{\bu\in \bbR^{|J|}\,:\,\|\bu\|_2 = 1}\frac{1}{\sqrt{n}} \|\bX_J \, \bu\|_2.
\eeq
Equivalently, $\lmin_{s}$ is the smallest singular value of among submatrices of $\frac{1}{\sqrt{n}}\bX$ made of at most $s$ columns.  Note that, indeed, $\Id(s)$ is equivalent to $\lmin_{s}>0$.

\subsubsection{Support recovery}

We now state the main result concerning the support recovery problem.  It states that, under $\Id((2+\es)\strue)$, the posterior distribution $\Pi$ concentrates sharply on the support of $\btrue$ --- which we assumed to be $\strue$-sparse --- as long as $\lambda$ and the nonzero coefficients are sufficiently large.

\begin{thm}\label{thm:support}
Consider a design matrix $\bX$, with $p \ge n$ and normalized column vectors \eqref{normalized}, that satisfies Condition $\Id((2+\es)\strue)$ for some fixed $\es>0$.
Assume that \eqref{lambda} holds and
\beq \label{rho}
\min_{j\in \Jtrue}|\beta_{\star,j}| \ge \rho := \frac{3 \sigma \sqrt{\lambda/n}}{\lmin_{(2+\es)\strue}} .
\eeq
Then, with probability at least $1-2p^{-c}$, $\Pi(\Jtrue) > \Pi(J)$ for all $J$, and in fact
\begin{align*}
\Pi\left(\Jtrue\right) \geq 1 - 4p^{-c}.
\end{align*}

\end{thm}

Under the conditions of \thmref{support}, some straightforward calculations imply that $\Jmap = \Jtrue$ with probability at least $1- 6p^{-c}$.  In particular, as $p \to \infty$, the MAP consistently recovers the support of the coefficient vector.  Note that the same is true for a subset drawn from~$\Pi$.
%In practice, direct simulation of the posterior $\Pi$ is in general not possible and we need to resort to MCMC techniques in order to obtain approximate realizations of $\Pi$. The practical aspect is developed in Section \ref{sec:num}.

%This support concentration result

The result applies in the ultra-high dimensional setting where $p$ is exponential in $n$, as long as the conditions are met.  Characterizing design matrices $\bX$ that satisfy $\Id((2+\es) \strue)$ in the ultra-high dimensional setting is an interesting open question beyond the scope of this paper.

We mention that, if $\strue$ is known and we restrict the prior over subsets $J$ of size exactly $\strue$, then the same conclusions are valid with $\es = 0$ and $\lmin_{(2+\es)\strue}$ replaced by $\lmin_{2\strue}$ in \eqref{rho}, yielding consistent support recovery under the minimum identifiability condition $\Id(2\strue)$.  In \secref{comparison}, we show that the Lasso estimator requires much more restrictive conditions on the design matrix and $\btrue$ to ensure it selects the correct variables with high probability.

Finally, we note that the concentration is even stronger.  Under the same conditions, if
\[\lambda = \frac{(1+\es) (23+5c) + m}{\es} \log p,\]
then
\beq \label{support-refine}
\sum_{J\subset [p]\,:\,  J\neq \Jtrue } |J|^m \Pi\left(J\right) \le 4p^{-c} \Pi(\Jtrue).
\eeq
We will use this refinement in the proof of \thmref{estimation-2}.

\subsubsection{Estimation}

Armed with results for the support recovery and prediction problems, we establish corresponding bounds for the estimation problem.  Our first result is a simple consequence of \thmref{prediction} and \prpref{support}.
%
%Define
%\beq \label{theta}
%\theta_{s} = \min_{J \subset [p]\,:\, |J|\leq s} \, \min_{\bu\in \bbR^{|J|}\,:\,\|\bu\|_\infty = 1}\frac{1}{\sqrt{n}} \|\bX_J \, \bu\|_2.
%\eeq
%When \thmref{prediction} and \prpref{support}, we easily deduce the following.
\begin{thm} \label{thm:estimation}
Consider a design matrix $\bX$ with $p \ge n$ and normalized column vectors \eqref{normalized}.
Assume $\lambda$ satisfies \eqref{lambda} with $\eps \le 1/2$.  Then with probability at least $1 - 3 p^{-c}$, we have
\[
\|\bbmap - \btrue\|_2 \le \sigma \sqrt{\frac{8 \strue \lambda}{n \lmin_{(2+\es)\strue}^2}}.
%\qquad \text{ and } \qquad
%\|\bbmap - \btrue\|_\infty \le \frac{{\rm PE}_{\rm map}}{\theta_{(2+\es)\strue}}.
\]
\end{thm}

We continue with bounds on the estimation error, this time in terms of the $l_{\infty}$-norm.
Based on \thmref{support} (and its proof), we deduce the following.
\begin{thm}\label{thm:estimation-1}
Let the conditions of \thmref{support} be satisfied.  Then, with probability at least $1 - 7p^{-c}$, we have
\begin{align}\label{est-ineq-1}
\|\bbmap - \btrue\|_\infty \le \sigma \sqrt{\frac{2 (c+1) \log p}{n\lmin^2_{\strue}}}.
\end{align}
\end{thm}
We emphasize that this estimator requires only the near minimum condition $\Id((2+\es)\strue)$ and that the nonzero components of $\btrue$ are somewhat larger than the noise level in \eqref{rho} to achieve the optimal (up to logs) dependence on $n,p$ of the $l_{\infty}$-norm estimation bound. We will develop this point further in our comparison with the Lasso.

We now study the performances of the posterior mean $\bbpm$ and that of the following variant
\begin{equation}
\widetilde \bbeta  = \sum_{J\subset [p]\,:\, \nu_J>0} \Pi(J) \bbhat_J, \qquad \nu_J := \min_{\bu\in \bbR^{|J|}\,:\, \|\bu\|_2 = 1} \frac{1}{\sqrt{n}}\|\bX_{J} \bu\|_2.
\end{equation}

Define the quantity $\nu_{\min} = \min_{J\subset [p]\,:\, \nu_J>0}\nu_{J}$, and note that $\nu_{\min}>0$.
%Set $\lambda_1 =\frac{(1+\es)(23+5c)+1 }{\es} \log p $ and $\lambda_2 = (62 + 4c) \log p$.

\begin{thm}\label{thm:estimation-2}
Let the conditions of \thmref{support} be satisfied and let $c\ge 1$. \\
\benum
\item Take $\lambda = \frac{(1+\es)(23+5c)+1 }{\es} \log p$.  Then, with probability at least $1-4p^{-c}$,
\begin{eqnarray*}
\|\widetilde \bbeta - \btrue\|_\infty &\le&  \sigma\sqrt{\frac{2(c+1) \log p }{ n\lmin^2_{\strue}}} + \frac{3}{\nu_{\min} \, p^c} \left[  \sigma\sqrt{(20+4c)\frac{\log p}{n}}  + \frac{\|\bX\btrue\|_2}{\sqrt{n}} + \nu_{\min}\|\btrue\|_\infty \right].
\end{eqnarray*}
\item If in addition $\Id(\smax)$ is satisfied,
\begin{eqnarray*}
\|\bbpm - \btrue\|_\infty &\le& \sigma\sqrt{\frac{2(c+1) \log p }{ n\lmin^2_{\strue}}} + \frac{3}{\nu_{\smax} \, p^c} \left[  \sigma\sqrt{(20+4c)\frac{\log p}{n}}  + \frac{\|\bX\btrue\|_2}{\sqrt{n}} + \nu_{\smax}\|\btrue\|_\infty \right].
\end{eqnarray*}
\item If in addition $\Id(\strue + \smax)$ is satisfied and $\lambda \ge (62 + 4c) \log p$,
\begin{align}\label{est-ineq-3}
\|\bbpm - \btrue\|_\infty \le \sigma\sqrt{\frac{2(c+1) \log p }{ n\lmin^2_{\strue}}} + \frac{2\sqrt{10} \sigma}{\sqrt{n}\nu_{\strue+\smax}} \left[  \frac{2\sqrt{\strue}}{p^c}  + \frac{1}{p^{\strue}} \right].
\end{align}
\eenum
\end{thm}

We note that $\bbpm$ requires at least $\Id(\smax)$.  (Recall that we assume $\smax$ is known such that $\strue \le \smax$.) In practice, when the sparsity is unknown, we make a conservative choice $\smax \gg 2\strue$ so that $\Id(\smax)$ is substantially more restrictive that $\Id(2\strue)$. Typically, we assume that $\strue = O\left(\frac{n}{\log p}\right)$ and we take $\smax$ of this order of magnitude. We will see below in \secref{gaussian} that for Gaussian design, the condition $\Id(\strue+\smax)$ is satified with probability close to $1$. On the other hand, the estimation result for $\widetilde \bbeta$ holds true under the near minimum condition $\Id((2+\es)\strue)$. For both estimators, their estimation bounds depend on the quantities $\nu_{\min}$, $\nu_{\smax}$, $\|\bX \btrue\|_2$ and $\|\btrue\|_\infty$ which can potentially yield a sub-optimal rate of estimation. Note however the presence of the factor $p^{-c}$ in the bound. In particular, if the nonzero components of $\btrue$ are sufficiently large, then the quantities $\nu_{\min}$, $\nu_{\smax}$, $\|\bX \btrue\|_2$ and $\|\btrue\|_\infty$ may be completely cancelled for a sufficiently large $c>0$. If $\Id(\strue+\smax)$ is satisfied, then we can derive a bound that no longer depends on $\|\bX \btrue\|_2$ and $\|\btrue\|_\infty$. We will also see below that this bound yields the optimal rate of $\l_\infty$-norm estimation (up to logs) for the estimator $\bbpm$ when the design matrix is Gaussian. Optimality considerations are further discussed in \secref{discussion} based on recent information bounds obtained elsewhere.

\subsubsection{Example: Gaussian design} \label{sec:gaussian}
%\noindent {\bf Example: Gaussian design.}
The quintessential example is that of a random Gaussian design, where the {\em row} vectors of $\bX$, denoted $\bx_1, \dots, \bx_n$, are independent Gaussian vectors in $\bbR^p$ with zero mean and $p\times p$ covariance matrix $ \bSigma$.  If we assume that $\bSigma$ has 1's on the diagonal, the resulting (random) design is just slightly outside our setting, since the columns vectors are not strictly normalized.  Our results apply nevertheless.  Therefore, it is of interest to lower-bound $\lmin_{s}$ for such a design.

We start by relating $\bX$ and $\bSigma$.  Consider $J\subset [p]$% with $|J| =s$
, and let $\bSigma_J$ denote the principal submatrix of $\bSigma$ indexed by $J$.  By \citep[Cor.~1.50 and Rem.~1.51]{vershynin}, there is a numeric constant $C>0$ such that, when $n \ge C |J|/\eta^2$, with probability at least $1- 2 \exp(-\eta^2 n/C)$, we have
\[
\left\|\frac1n \bX_J^\top \bX_J - \bSigma_J\right\| \le \eta \|\bSigma_J\|,
\]
%\begin{thm}[\citep[Cor.~1.50 and Rem.~1.51]{vershynin}] \label{thm:vershynin}
%There is a numeric constant $C > 0$ such that, when $\bx_1, \dots, \bx_n \sim \cN(
%\end{thm}
where $\|\cdot\|$ denotes the matrix spectral norm. When this is the case, by Weyl's theorem \citep[Cor.~IV.4.9]{MR1061154},
\[
\lambda_{\rm min} \left(\frac1n \bX_J^\top \bX_J\right) \ge \lambda_{\rm min} (\bSigma_J) - \eta \lambda_{\rm max}(\bSigma_J),
\]
where $\lambda_{\rm min}(\bA)$ and $\lambda_{\rm max}(\bA)$ denote the smallest and largest eigenvalues of a symmetric matrix $\bA$.  Define
\[
\eta_\bSigma(s) = \max_{J : |J| \le s} \frac{\lambda_{\rm max}(\bSigma_J)}{\lambda_{\rm min} (\bSigma_J)}, \qquad \lambda_\bSigma(s) = \min_{J : |J| \le s} \lambda_{\rm min} (\bSigma_J).
\]
Assume that
\[
n \ge \frac{a C s \log p}{\eta_\bSigma(s)^{2}},
\]
for some $a \ge 2$.  Then, with probability at least $1 - 2 p^{-a/2}$,
\[
\lmin_{s} \ge \frac{\lambda_\bSigma(s)^{1/2}}{2}.
\]
For example, in standard compressive sensing where $\bSigma$ is the identity matrix, we have $\eta_\bSigma(s) = \lambda_\bSigma(s) = 1$ for all $s$, in which case with high probability $\lmin_{s} \ge 1/2$ when $n \ge 2 C s \log p$. Consequently, the $l_\infty$-norm estimation bounds in \eqref{est-ineq-1} and \eqref{est-ineq-3} are of the order $b \sigma \sqrt{\log(p)/n}$ for some numerical constant $b>0$. Again, the constants are loose in this discussion.

\subsection{Bayesian variable selection with an independence prior} \label{sec:bayes}
Many Bayesian techniques for model selection have proposed in the literature; see \citep{MR2000752} for a comprehensive review.  That same paper suggests a procedure similar to ours, except that it is a bonafide Bayesian model and they use the following independence sparsity prior
\[
\tilde{\pi}(J) = \omega^{|J|} (1 - \omega)^{p - |J|},
\]
where $\omega \in (0,1)$ controls the sparsity level.  Roughly, $\lambda$ for our prior corresponds to $\log(1-1/\omega)$ for this prior.  It so happens that, the same arguments lead to the same results.
%, meaning that Theorems \ref{thm:support} and \ref{thm:estimation} are valid under this prior.
Also, as argued in \citep[Sec.~3.3]{MR2000752}, the fully-specified Bayesian model with prior $\bbeta_J \sim \cN({\bf 0}, (\bX_J^\top \bX_J)^{-1})$ is very closely related to our exponential weights method.

\subsection{$\ell_0$-penalized variable selection} \label{sec:bic}
\cite{MR2443189} not only showed that BIC was consistent when $p < \sqrt{n}$ (under some mild conditions on the design matrix), they also suggested a modification of the penalty term to yield a method that is consistent for larger values of $p$ when the number of variables in the true (i.e., sparsest) model $\strue$ is bounded independently of $n$ or $p$.

By a simple modification of our arguments, our results for the exponential weights MAP is seen to apply to
\[
\Jhat = \argmin_{J: |J| \le \smax} \ \by^\top (\bI -\bP_J) \by + \lambda |J|.
\]
Consequently, our work extends that of \cite{MR2443189} to the case where $\strue$ increases with~$p$.

\section{Comparison with the Lasso and {\sc mc+}}
\label{sec:comparison}

In this section, we compare the theorethical performances of our
procedure with other well-known $l_\infty$-estimation and support
recovery techniques used in high-dimensional variable selection.

\subsection{Lasso}

The Lasso estimator is the solution of the convex minimization problem
$$
\bbhat^{L} = \argmin_{\bbeta\in\bbR^p} \left\lbrace \frac{1}{n}\|\by
- \bX\bbeta\|_2^2 + 2\lambda_L \|\bbeta\|_1  \right\rbrace,
$$
where $\lambda_L = A \sigma \sqrt{\log(p)/n}$, $A>0$ and $\|\cdot\|_1$ is the $l_1$-norm. The Lasso has
received considerable attention in the literature over the last few
years \citep{BoBach08,B07,BTW07,MB06,MY06,ZY06}. It is not our goal
to make here an exhaustive presentation of all existing results. We
refer to Chapter 4 in \citep{louthese} and the references cited
therein for a comprehensive overview of the literature. %On the theorethical side,
%sparsity oracle inequalities have been established for the Lasso
%under
% a variety of assumptions on the design matrix $\bX$. Concerning the prediction and $l_q$-norm estimation problems with $1\leq q<\infty$, one widely used condition on $\bX$ is
% the Restricted Isometry Condition (or some related version). See for instance \citep{CT07,BRT07}.
%%assume that $\bX$ satisfies a restricted isometry condition $\mathbf{RI}(\strue,c_0)$ (or some related versions)

Concerning the $l_\infty$-norm estimation and support recovery
problems, the most popular assumption is the Irrepresentable
Condition \citep{BoBach08,louthese,MY06,W06,ZY06} denoted from now
on by $\textbf{IC}(\strue)$. See for instance Assumption 4.2 in
\citep{louthese}. The condition $\textbf{IC}(\strue)$ is strictly
more restrictive than the identifiability $\Id(2\strue)$ and does
not hold true in general when the columns of the design matrix $\bX$ are not weakly
correlated. Define $d_\star = \| \Psi_\star^{-1}\vsign(\btrue)\|_\infty$ where
$\Psi_\star:=\frac{1}{n}\bX_{\Jtrue}^{\top}\bX_{\Jtrue}$.
%and $\overline\btrue \in \bbR^{\strue}$ is obtained by keeping only the nonzero components of $\btrue$.
The following result is the key to our analysis. Let
$\textbf{IC}(\strue)$ holds true and let the nonzero components of
$\btrue$ be sufficiently large: $\min_{j\in
\Jtrue}|\beta_{\star, j}| > A\sigma d_\star \sqrt{\log(p)/n}$.
 Then, with probability at least
$1-2p^{1-\frac{ A^{2}}{16}} - \strue p^{-\frac{A^{2}}{2}}$, the
Lasso solution is unique and satisfies
\begin{equation}\label{theo-sup-norm-chap-3}
c \sigma d_\star \sqrt{\frac{\log p}{n}} \leq \|\bbhat^{L} - \btrue\|_\infty \leq C\sigma d_\star \sqrt{\frac{\log p}{n}},
\end{equation}
for some numerical constants $C\geq c>0$ that can depend only on $A$. See Theorem 4.1 in
\citep{louthese} for a more precise statement.

We say that a $l_\infty$-norm estimation rate is optimal if it is of
the form $\alpha\sigma \sqrt{\log(p)/n}$ where $\alpha>0$ is
an absolute constant as in the case of gaussian sequence model
($n=p$ and $\bX = \bI_n$ the $n\times n$ identity matrix).
In view of the previous display, the Lasso does not attain in most
cases the optimal $l_\infty$-norm estimation rate. Indeed, the
quantity $d_\star$ generally depends on $\strue$ unless the
correlations between the columns of the design matrix $\bX$ are very weak. Consider for the instance the case where $\min_{j\neq k}|(\Psi_{\star}^{-1})_{j,k}|\ge \rho$ for some fixed $\rho>0$. Then, we can easily find a $\strue$-sparse vector $\btrue$ such that $d_\star\geq \rho \strue$ and the $l_\infty$-norm estimation rate of the Lasso is then suboptimal by a factor $\strue$.

% The exception is when the design matrix $\bX$ admits weakly correlated columns. Assume for instance that $\Psi$ satisfies a mutual coherence condition
% $\max_{i\neq j}|\Psi_{i,j}|\leqslant \frac{1}{cs}$ for some $c>1$. Then, elementary computations give that $d^{*} \leqslant \frac{c}{c-1}$ and consequently
% the Lasso achieves the optimal $l_\infty$-norm estimation rate in this particular situation provided that $\rho>Cr$ where $C>0$ is a constant independent of $s$.
%See Section 4.3 in \citep{louthese} for more details.

Unlike the Lasso, our exponential weights procedure does not suffer
from this limitation. Indeed, our procedure achieves the optimal
$l_\infty$-norm estimation rate and support recovery provided that
Condition $\Id((2+\es)\strue)$ holds true, which can be the case even for design matrices $\bX$ with strongly correlated columns. %and that the nonzero components of $\btrue$ are larger than the noise level.

\medskip
\noindent {\bf Gaussian design.}
Consider the Gaussian design of \secref{gaussian}, but assume now that $\bSigma = I_{p\times p}$ the $p\times p$ identity matrix.
Although the design $\frac{1}{\sqrt{n}}\bX$ satisfies the restricted
isometry with probability close to $1$, there is no guaranteed that
$\bX$ also satisfies an irrepresentable condition
$\mathbf{IC}(\strue)$. Let's assume that this is the case
for the sake of comparison. Then, we can show
with probability close to $1$ that $\Psi_\star^{-1}$ satisfies the mutual
coherence condition $\max_{j\neq k} |(\Psi_\star)^{-1}_{j,k}|\leq
\frac{1}{\sqrt{\strue}}$ where the dependence on $\strue$ cannot be
improved. Thus, we get $d_\star \le \sqrt{\strue}$ and we cannot guarantee the optimality of the $l_\infty$-norm estimation bound for the Lasso under the irrepresentable condition. Consequently, we
need the condition $\min_{j\in
\Jtrue}|\beta_{\star, j}| \ge C\sigma\sqrt{\strue \log(p)/n}$ for some absolute constant $C>0$ in order to guarantee exact support recovery for the Lasso. This condition is to be compared to \eqref{rho} for the exponential weights estimators. In that case, we have $\nu_{\tilde s}>1/2$ with probability close to $1$ when $\tilde s = O\left(n/\log p\right)$, so that \eqref{rho} becomes simply $\min_{j\in
\Jtrue}|\beta_{\star, j}| \ge C \sigma \sqrt{\log(p)/n}$ for some numerical constant $C>0$. This condition is less restrictive than that for the Lasso by a factor $\sqrt{\strue}$. Next, we note also that for a Gaussian design, the estimation bounds \eqref{est-ineq-1} and \eqref{est-ineq-3} for the exponential weights estimators are optimal (up to log) whereas the estimation bound for the Lasso contains the additional factor $\sqrt{\strue}$.

\medskip

Recently, in the framework of instrumental regression, \cite{GauTsy11} established for an $l_1$-norm minimization procedure $\hat\bbeta^{D}$ close to the Dantzig selector (see (3.5) there) that with probability close to $1$
$$
\|D_{\bX}^{-1}(\hat \bbeta^{D} - \btrue)\|_q \leq 2\sigma\sqrt{\frac{\log p }{\kappa^2_{q,\Jtrue} n}},
$$
where the sensitivity
$$\kappa_{q,\Jtrue}:=\inf_{\Delta \in C_{\Jtrue}:\|\Delta\|_q=1}\left| \frac{1}{n}D_{\bX}\bX^\top \bX D_{\bX}\Delta\right|,$$
with $C_{\Jtrue} = \left\lbrace \Delta\in \bbR^p\,:\, \| \Delta_{\Jtrue^c} \|_1 \leq \frac{1+c}{1-c}\|\Delta_{\Jtrue}\|_1 \right\rbrace$ for some $0<c<1$, $D_{\bX} = \mathrm{diag}(\bX_{1*},\cdots,\bX_{p*})$, $\bX_{k*} = \max_{1\leq i \leq n}|\bX_k^{(i)}|$ for any $1\leq k \leq p$ and the $\bX_k^{(i)}$ are the components of $\bX_k$. An enticing property of the sensitivity approach is that the quantities $\kappa_{q,J}$ can be computed in reasonable time for small $J$, yielding a computationally tractable procedure to build confidence interval
for the estimation of $\btrue$. The downside is that without any further conditions on $\bX$, the optimal dependence of the bound on $\strue$ is not clear.
For instance, assume in addition that $\bX$ satisfies a restricted eigenvalue condition as in \citep{BRT07}, then the dependence of the above bound on $\strue$ can be proved
 to be optimal for any $1\leq q\leq 2$ (See Section 9 in \cite{GauTsy11}). However, if we only assume that the condition
$\Id((2+\eps)\strue)$ is satisfied, then we cannot establish a clear
comparison between the exponential weights and $\hat\bbeta^D$. The
exponential weights estimator achieves the optimal estimation rate
while the dependence on $\strue$ of the $l_q$-norm estimation bound
is not explicit for $\hat\bbeta^D$.

\subsection{{\sc mc+}}

The {\sc mc+} estimator initially proposed by \cite{MR2604701} is the solution of the following nonconvex minimization problem:
\begin{equation}
\bbhat^{{\sc MC+}} = \mathrm{argmin}_{\bbeta \in \bbR^p} \left\lbrace \frac{1}{n}\|\by-\bX\bbeta\|_2^2 + \sum_{j=1}^p \Upsilon(|\bbeta_j|,\lambda_{\sc MC},\gamma)   \right\rbrace,
\end{equation}
where $\lambda_{MC},\gamma>0$ and the {\sc mc+} penalty function $\Upsilon$ is nonconvex, equal to $0$ outside a compact neighborhood of $0$ and admits a nonzero right derivative at $0$.  See equations (2.1)-(2.3) in \citep{MR2604701} for more details.

The performance of this estimator is established in Theorem 1 of \citep{MR2604701}, where the tuning of the parameter $\lambda_{\sc MC}$ requires the knowledge of $\strue$ (which is $d^o$ in that paper) and the optimal theoretical choice of $\gamma$ is proportional to $\nu_{\bar s}^{-1}$ (which is $d^*$ in that paper).  Let us also emphasize that the
choice $\gamma \propto \nu_{\bar s}^{-1}$ requires that $\Id(\bar s)$ is satisfied where $\bar s$ is in practice a conservative upper bound on $\strue$. In addition, this quantity $\nu_{\bar s}$ is delicate to compute in practice. Note that the exponential weights do not present the same limitations. Indeed, no prior knowledge of $\strue$ is required and we only need the condition $\Id((2+\es)\strue)$ for an arbitrarily small $\es>0$ to establish the consistency of $\bbmap$ even if the parameter $\bar s$ is chosen conservatively (for example, $\bar s = [n/2]$ if no other information is available). In addition, the tuning of the parameters for the exponential weights do not require to compute any restricted eigenvalues.  We mention that the assumptions in Theorem 1 of \citep{MR2604701} do not guarantee the identifiability of $\btrue$.

%
% \subsection{MCP and SCAD penalty}
%
% Sparse Riesz condition
%
% \subsection{Greedy}
%
%
% \subsection{OMP}
%
% Forward stepwise selection, or orthogonal matching pursuit (OMP), is too short-sighted and is easily mislead if there are substantial correlations in the predictors.  For a simple example, consider the case where $\bX_1, \dots, \bX_{p-1}$ are orthogonal while
% \[\bX_p = (1-\eta) (\bX_1 + \bX_2) + \sqrt{\frac{\eta}{p-3}} (\bX_3 + \cdots + \bX_{p-1})\]
% for some $\eta \in (0,1)$, $\Jtrue = \{1,2\}$ with $\beta_1^* = \beta_2^*$, and we take $\es < 1$.  Note that the conditions in Section XXX are satisfied.  Suppose there is no noise.  In the first step, OMP chooses the predictor that is the most correlated with the response, which here is $\beta_1^* (\bX_1 + \bX_2)$, which is $\bX_p$ when $\eta < 1 - 1/\sqrt{2}$, so the model chosen by OMP includes $\bX_p$.
%
% \eac{Le probleme que je vois c'est que la routine Metropolis de notre algo se trompe aussi, il me %semble...}

\section{Numerical Experiments}
\label{sec:num}

In this section, we illustrate the performance of
the procedure (\ref{posterior-sparse}) in variable selection and
$l_\infty$-norm estimation on a simulated data set. The posterior \eqref{posterior-sparse} is simulated via MCMC. In a
nutshell, we construct an ergodic Markov chain $(\bbeta_t)_{t\ge 0}$ with invariant probability distribution the posterior
(\ref{posterior-sparse}). Then, we get from \citep{RobCas04} that
$$
\lim_{T\rightarrow \infty} \frac{1}{T}\sum_{t=T_0+1}^{T_0 + T}
\bbeta_t = \bbpm,\quad
\pi-a.s.\,,
$$
where $T_0\geq 0$ is an arbitrary number. In practice, we use
$T_0=3000$ and $T=7000$.
% and
%$$
%\lim_{T\rightarrow \infty} \frac{1}{T}\sum_{t=T_0+1}^{T_0 + T} \1 \{J_T = \Jtrue\} = \Pi(\Jtrue).
%$$

We refer to \citep{AL11,MR2816337} for more
details on the computational aspect. Our numerical study follows
those carried out in these references except that we concentrate on
the $l_\infty$-norm estimation and variable selection performances
of the procedure $\bbpm$. Note that the exponential weights procedures
considered in the present paper and in \citep{AL11,MR2816337} differ only through the  tuning of the parameters. \citep{AL11,MR2816337} consider indeed the prediction problem whereas we concentrate on the $l_\infty$-norm estimation and support recovery problems, which require a different tuning to guarantee the theoretical consistency. Note also that $\bbpm$ is not sparse since it is obtained as the expectation of the posterior $\Pi$. However, in view of \thmref{estimation-2}, a simple thresholding of $\bbpm$ with a threshold of the order of the noise level yields consistent support recovery. In our simulations, we observe that few components of $\bbpm$ are significantly far from $0$ whereas the remaining ones are extremely small thus making the choice of the threshold easy in practice. From now on, we will denote indifferently by {\sc aew} the procedure $\bbpm$ and the thresholded $\bbpm$.

Following the numerical experiments of \citep{MR2382644}, we
consider the model (\ref{model}) where $\bX$ is an $n\times p$
matrix with independent standard Gaussian entries, the target vector
$\btrue = \1_{\{j\leq \strue\}}$ for some fixed $\strue\geq 1$ and the
noise variance satisfies $\sigma^2 = \|\bX \btrue \|^2/(9n)$. For each different setting of $(n,p,\strue)$, we perform $100$
replications of the model and compare our estimator {\sc awe} with
other procedures in the literature on sparse estimation:
\begin{enumerate} \setlength{\itemsep}{0in}
  \item The Lasso estimator; % with regularization parameter obtained by ten-fold cross-validation;
  \item The {\sc mc+} estimator of  \citep{MR2604701}; % with regularization parameter obtained by ten-fold cross-validation;
  \item The {\sc scad} estimator of \citep{FanLi01}. % with regularization parameter obtained by ten-fold cross-validation.
\end{enumerate}

All these procedures are readily implemented in R. We use the glmnet
package to compute the Lasso estimor and the ncvreg package to
compute the {\sc scad} and {\sc mc+} estimators. The {\sc awe}
estimator was computed through the MCMC algorithm described in
Section 7 in \citep{MR2816337}.

\begin{figure}[h]\label{fig1}
\psfrag{a}[tr][][.8][45]{{\sc aew}}
\psfrag{b}[tr][][.8][45]{Lasso}
\psfrag{c}[tr][][.8][45]{{\sc mc+}}
\psfrag{d}[tr][][.8][45]{{\sc
scad}}
\psfrag{f}[][][.8][0]{$\|\hat \bbeta -\bbeta_*\|_\infty$}
\begin{center}
\includegraphics[width=0.5\textwidth]{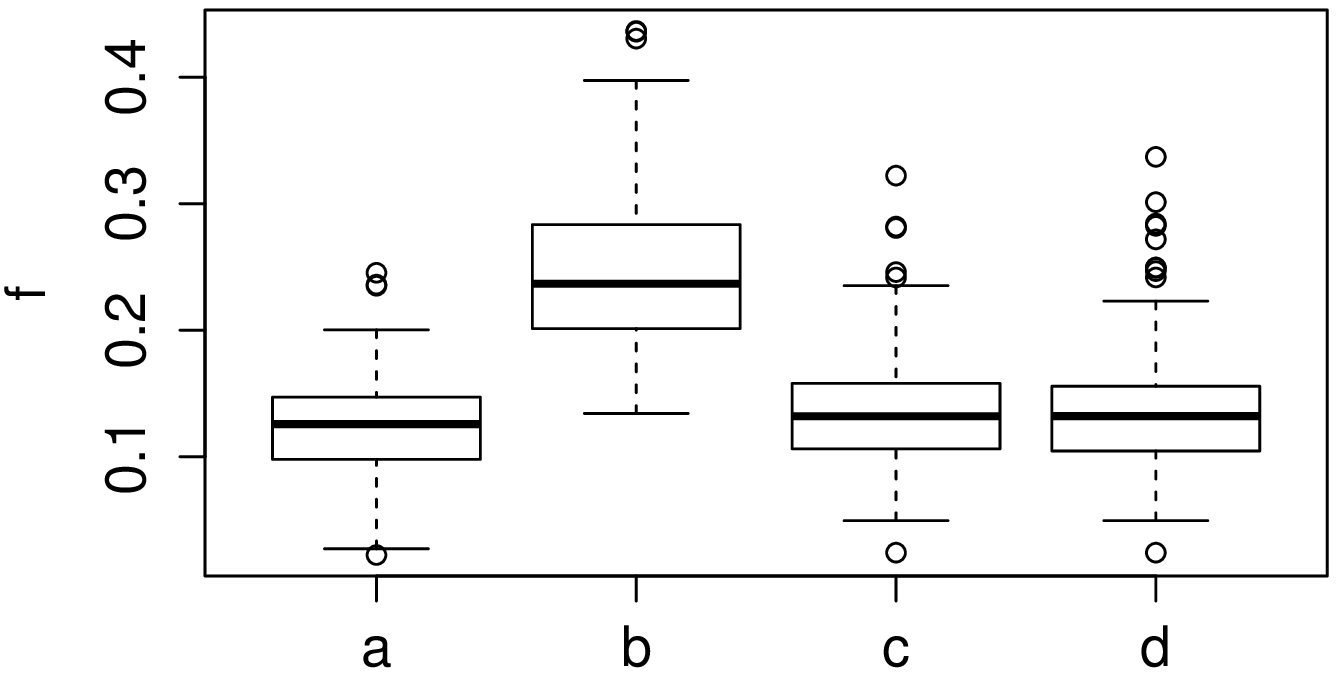}%
\includegraphics[width=0.5\textwidth]{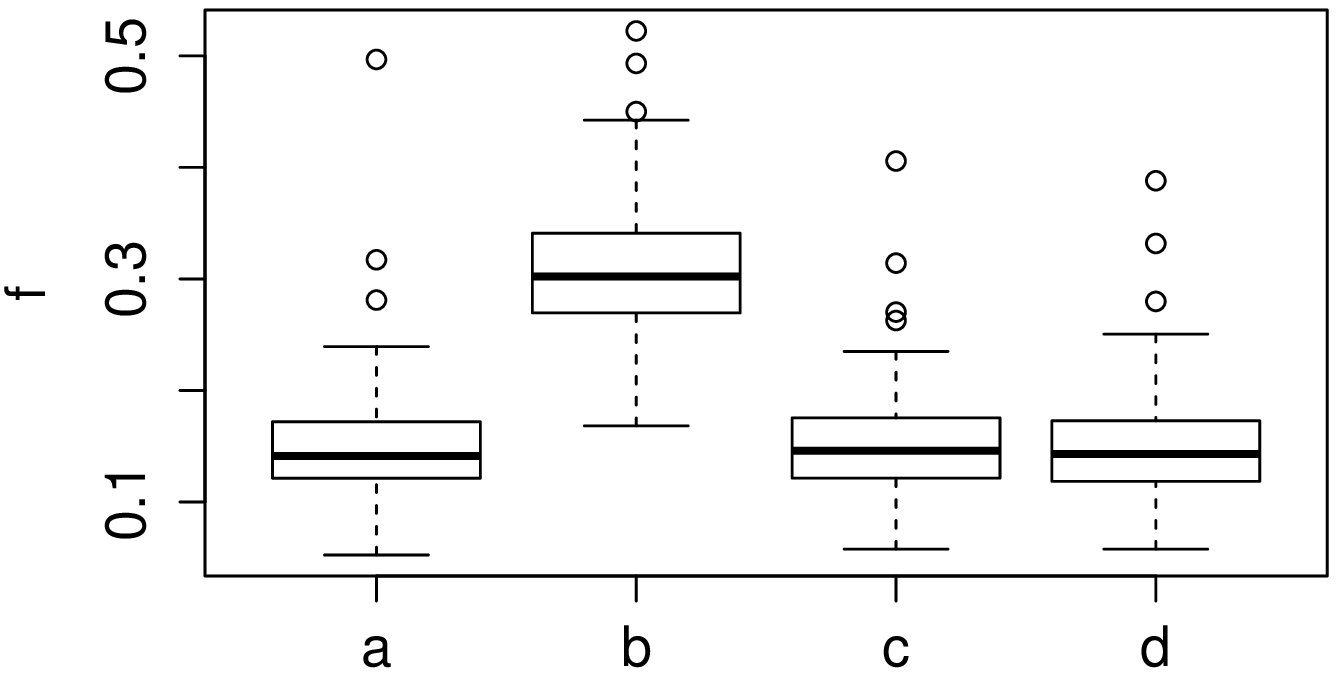}%
\end{center}
%\vspace{-1.5cm}
%\begin{center}
%\includegraphics[width=0.5\textwidth]{boxplotpred12.eps}%
%\includegraphics[width=0.5\textwidth]{boxplotpred25.eps}%
%\end{center}
\caption{Independent Gaussian design. Boxplots of estimation performance measure $\|\hat \bbeta
-\btrue\|_\infty$ over 100 realizations for the {\sc
aew}, Lasso, {\sc mc+} and {\sc scad} estimators. {\it Left:}
$(n,p,\strue)=(100,200,5)$. {\it Right:}
$(n,p,\strue)=(200,1000,10)$.
} \label{FIG:boxplots}
\end{figure}

Figures \ref{FIG:boxplots} contains the comparative boxplots for the $l_\infty$-norm
error over the $100$ repetitions for the independent Gaussian design. Table \ref{TAB:tabest} contains the average $l_\infty$-norm error and the standard deviation over the $100$ repetitions for the independent Gaussian design. We observe that the {\sc awe} estimator outperforms the Lasso estimator and exhibit performances similar to {\sc mc+} and {\sc scad}.

%{\it Bottom:} Prediction performance: $|\bX(\hat\theta -\theta^*)|_2^2/n$.

\begin{table}[h]
\begin{center}
\begin{tabular}{c|l|l|l|l}
$(n,p,\strue)$ &  {\sc aew} & Lasso & {\sc mc+} & {\sc
scad}\\
\hline \hline $(100,200, 5)$&  \textbf{ \footnotesize 0.124} &{ \footnotesize 0.249} & {\footnotesize 0.137} &{ \footnotesize 0.138}   \\
&   {\textbf{\footnotesize (0.041)}}  &  {\footnotesize (0.068)  } & {\footnotesize (0.050) }  & {\footnotesize (0.056) } \\
\hline $(200,1000, 10)$ &  {\footnotesize
0.151}  &  {\footnotesize 0.309  } & {\footnotesize 0.153 }  &
{\textbf{\footnotesize 0.149}}
\\
&     {\footnotesize (0.055)}  &  {\footnotesize (0.063)  } & {\footnotesize (0.051) }  &  {\textbf{\footnotesize (0.050)} }\\
\end{tabular}
\end{center}
\caption{Independent Gaussian design. Means and standard deviations of performance measures over
100 realizations for the {\sc aew}, Lasso, {\sc mc+} and
{\sc scad} estimators.} \label{TAB:tabest}
\end{table}

We note in our simulation study that the four procedures always
select the $\strue$ active covariates but also select non-active
ones. Table \ref{TAB:varsel} contains the average support recovery
false positive rate over the $100$ repetitions for the four
procedures considered in this study. We observe that the Lasso tends
to select too many covariates as was already known. The
(thresholded) {\sc awe} estimator outperforms all other procedures
in the support recovery problem.

\begin{table}[h]
\begin{center}
\begin{tabular}{c|l|l|l|l}
$(n,p,\strue)$ &  {\sc awe} & Lasso & {\sc mc+} & {\sc
scad}\\
\hline \hline $(100,200, 5)$
&     {\textbf{\footnotesize 1.60 }}  &  {\footnotesize 21.55  } & {\footnotesize 1.75 }  & {\footnotesize 3.02 }  \\
\hline $(200,1000, 10)$ &    {\textbf{\footnotesize 1.98}}  &
{\footnotesize 51.88 } & {\footnotesize 2.49 }  &  {\footnotesize
5.22
 }\\
\end{tabular}
\end{center}
\caption{Average support recovery
false positive rate over
100 realizations for the {\sc aew}, Lasso, {\sc mc+} and
{\sc scad} estimators.}\label{TAB:varsel}
\end{table}

\section{Discussion}
\label{sec:discussion}

We established some performance bounds for exponential weights when
applied to solving the problems of {\em prediction}, {\em
estimation} and {\em support recovery}, and deduced similar results
for a slightly different Bayesian model selection procedure
\citep{MR2000752} and $\ell_0$-penalized (BIC-type) variable
selection.  How sharp are these bounds?  We did not optimize the
numerical constants appearing in our results, simply because we
believe our bounds are loose and also because there are no known
sharp information bounds for theses problems, except in specific
cases \citep{jin2012optimality}.  That said, there are some results
available in the literature
\citep{verzelen,raskutti2009minimax,LPTV} and our bounds come close
to these.  For example, from \citep{raskutti2009minimax} we learn
that, when $\Id(2\strue)$ holds, there is a universal constant $C >
0$ such that, for any estimator $\bbhat$ that knows $\strue$,
\[
\|\bbhat - \btrue\|_2 \ge C \sigma \sqrt{\frac{\strue \log (p/\strue)}{n \kappa_{2\strue}^2}}
\]
with probability at least $1/2$, where
\beq \label{lmax}
\kappa_s := \max_{J \subset [p]\,:\, |J|\leq s} \, \min_{\|\bu\| = 1} \frac1{\sqrt{n}} \|\bX_J \, \bu\|_2;
\eeq
and from \citep{verzelen}, we learn that, for another universal constant $C'>0$,
\[
\E \|\bbhat - \btrue\|_2^2 \ge C' \sigma^2 \left(\frac{\strue \log (e p/\strue)}{\kappa_{2\strue}^2} \vee \frac1{\lmin_{2 \strue}^2}\right).
\]
Thus we see that our estimation bounds \eqref{est-ineq-1} and \eqref{est-ineq-3} come quite close to these information bounds.  Of course, there is a trade-off with computational tractability, as computing the exponential weights estimates (of even approximating them) in polynomial time remains an open problem.  That said, the numerical experiments show that these methods are promising.

\section{Proofs} \label{sec:proofs}
For the sake of brevity, we let $\|\cdot\| = \|\cdot\|_2$ throughout this section.

\subsection{Proof of \thmref{prediction}}

Define $\bxi_J = \bP_J(\by) -\bX \btrue$.  For $J \subset [p]$ with $|J| = s$, we have
\beq \label{ratio}
\frac{\Pi(J)}{\Pi(\Jtrue)}
= \frac{\binom{p}{\strue}}{\binom{p}{s}} \exp\left(\lambda(\strue-s) +\frac{1}{2 \sigma^2}(\|\bP_{\Jtrue}^{\perp}(\bz)\|^2 -\|\bP_J^{\perp}(\by)\|^2)\right)
\eeq
with
\beq \label{quad1}
\|\bP_{\Jtrue}^{\perp}(\bz)\|^2 -\|\bP_J^{\perp}(\by)\|^2
= 2 \bz^T (\bxi_J - \bxi_{\Jtrue}) + \|\bxi_{\Jtrue}\|^2 -\|\bxi_{J}\|^2 .
\eeq
For the inner product on the RHS, note that $\bxi_J \in {\rm span}(\bX_{J \cup \Jtrue})$ and $\bxi_{\Jtrue} \in {\rm span}(\bX_{\Jtrue})$, so that
\beq \label{inner1}
\big|2 \bz^T (\bxi_J - \bxi_{\Jtrue})\big| = \big| 2 (\bP_{J \cup \Jtrue} \bz)^T (\bxi_J - \bxi_{\Jtrue}) \big| \le 2 \|\bP_{J \cup \Jtrue} \bz\| \ \|\bxi_J - \bxi_{\Jtrue}\|,
\eeq
by Cauchy-Schwarz's inequality.

\begin{lem} \label{lem:proj}
For any $c>0$, with probability at least $1 - p^{-c}$,
\beq \label{proj}
\|\bP_{J}\, \bz\|^2 \le (20 + 4 c) \sigma^2 |J| \log p, \quad \forall J \subset [p].
\eeq
\end{lem}

Set $\zeta_J = \sqrt{(20 + 4 c) (|J|+\strue) \log p}$. Using \lemref{proj} in \eqref{inner1}, from \eqref{quad1} we have
\begin{eqnarray}
\|\bP_{\Jtrue}^{\perp}(\bz)\|^2 -\|\bP_J^{\perp}(\by)\|^2
&\le& \sigma \zeta_J \, \|\bxi_J - \bxi_{\Jtrue}\| + \|\bxi_{\Jtrue}\|^2 -\|\bxi_{J}\|^2 \notag \\
&\le& \sigma \zeta_J \, \big(\|\bxi_J\| + \|\bxi_{\Jtrue}\|\big) + \|\bxi_{\Jtrue}\|^2 - \|\bxi_{J}\|^2 \notag \\
&\le& 4 \sigma^2 \zeta_J^2 + \frac32 \|\bxi_{\Jtrue}\|^2 - \frac12 \|\bxi_{J}\|^2 \notag \\
&\le& 6 \sigma^2 \zeta_J^2 - \frac12 \|\bxi_{J}\|^2, \label{zeta}
\end{eqnarray}
where we used the identity $ab \le 2 a^2 + b^2/2$ in the third inequality, and \lemref{proj} to bound $\|\bxi_{\Jtrue}\|^2$ in the last inequality.

We tackle the first part.  By definition, $\Pi(\Jmap) \ge \Pi(\Jtrue)$.  Take any $J$ such that $\Pi(J) \ge \Pi(\Jtrue)$ and let $s = |J|$.  Plugging in the bound \eqref{zeta} into \eqref{ratio}, and using some crude bounds, we have
\beqn
1 \le \frac{\Pi(J)}{\Pi(\Jtrue)} &\le& \exp\left(\strue \log p + \lambda (\strue-s) + 3 (s+\strue) (20 + 4c) \log p - \frac1{4\sigma^2} \|\bxi_{J}\|^2 \right) \\
&\le& \exp\left(\strue \big(\lambda + (61 + 12c) \log p \big) - \frac1{4\sigma^2} \|\bxi_{J}\|^2 \right),
\eeqn
where we used the fact that $\lambda \ge (62 + 12c) \log p$ in the last inequality.
This in turn implies
\[
\|\bxi_{J}\|^2 \le 4\sigma^2 \cdot \big( \lambda \strue + (61 + 12c) \log p \big) \le 8 \sigma^2 \lambda,
\]
and the first part of \eqref{pred} follows from that.

We now turn to the second part.  Define $\cJ = \{J: \|\bxi_J\| > \sigma \sqrt{10 \strue \lambda}\}$.
We have
\begin{eqnarray}
\|\bX \bbpm -\bX \btrue\| &\le& \sum_J \|\bxi_J\| \Pi(J) \notag\\
&\le& \sigma \sqrt{10 \lambda \strue} \sum_{J \notin \cJ} \Pi(J)+ \sum_{J \in \cJ} \|\bxi_J\| \frac{\Pi(J)}{\Pi(\Jtrue)}.\label{thm-pred-interm-1}
\end{eqnarray}

By \eqref{ratio} and \eqref{zeta}, we have
\beqn
\|\bxi_J\| \ \frac{\Pi(J)}{\Pi(\Jtrue)}
&\le& \|\bxi_J\| \ \frac{\binom{p}{\strue}}{\binom{p}{s}} \exp\left(\lambda(\strue-s) + 3 \zeta_J^2 - \frac1{4\sigma^2} \|\bxi_{J}\|^2 \right) \\
&\le&  \frac{\sqrt{10} \sigma}{\binom{p}{s}} \exp\left(\lambda(\strue-s) + \strue \log p + 3 \zeta_J^2 - \frac1{5\sigma^2} \|\bxi_{J}\|^2 \right),
\eeqn
where we used the fact that $x e^{-x^2} \le 1/\sqrt{2}$ for all $x$, and $\binom{p}{\strue} \le p^{\strue}$.
Hence, since $\lambda \ge (62 + 4c) \log p$, we have
\begin{eqnarray}
\sum_{J \in \cJ} \|\bxi_J\| \frac{\Pi(J)}{\Pi(\Jtrue)}
&\le& \sum_{s=0}^{\smax} \sum_{J: |J|=s} \frac{\sqrt{10} \sigma}{\binom{p}{s}} \exp\left(\lambda(\strue-s) + \strue \log p + 3 \zeta_J^2 - 2 \lambda \strue \right) \notag \\
&\le& \sqrt{10} \sigma \sum_{s=0}^{\smax} \exp\left(- (\strue +s) (\lambda - (61 + 12c) \log p) \right) \notag\\
&=& \sqrt{10} \sigma \cdot 2 \exp\left(- \strue (\lambda - (61 + 12c) \log p) \right) \notag \\
&\le& 2 \sqrt{10} \sigma p^{-\strue}. \label{sum-cJ}
\end{eqnarray}
The result now follows from
\[
\sigma \sqrt{10 \lambda \strue} + 2 \sqrt{10} \sigma p^{-\strue} \le \sqrt{10} \sigma (\sqrt{\lambda\strue} + 1) \le \sigma \sqrt{12 \lambda\strue},
\]
since $p \ge 2$ and $\strue \ge 1$, as well as $\lambda \ge 25$.

\subsection{Proof of \prpref{support}}
Remember \eqref{ratio}.  We reformulate \eqref{quad1} in the following way
\begin{eqnarray}
\|\bP_{\Jtrue}^{\perp}(\bz)\|^2 -\|\bP_J^{\perp}(\by)\|^2
&=& \by^\top (\bP_J - \bP_{\Jtrue}) \by \notag \\
&=& - \|\bP_J^\perp \bX \btrue\|^2 - 2 \< \bP_J^\perp \bX \btrue, \bz \> + \bz^\top (\bP_J - \bP_{\Jtrue}) \bz. \label{fidel}
\end{eqnarray}
Let $\cJ_{s,t} = \{J \subset [p]: |J| = s, |J \cap \Jtrue| = t, J \neq \Jtrue\}$.
We first bound the inner product in \eqref{fidel}.

\bigskip

\begin{lem} \label{lem:inner}
For any $c>0$, with probability at least $1 - p^{-c}$,
\beq \label{inner}
\frac{\< \bP_J^\perp \bX \btrue, \bz \>^2}{\|\bP_J^\perp \bX \btrue\|^2} \le (10+2c) \sigma^2 (s \vee \strue - t) \log p,
\eeq
for all $J \in \cJ_{s,t}$ with $t \le s \wedge \strue$.
\end{lem}

\bigskip
We now bound the quadratic term in \eqref{fidel}.
\begin{lem} \label{lem:quad}
For any $c>0$, with probability at least $1 - p^{-c}$,
\beq \label{quad}
\bz^\top (\bP_J - \bP_{\Jtrue}) \bz \le (20+4c) \sigma^2 (s \vee \strue - t) \log p,
\eeq
for all $J \in \cJ_{s,t}$ with $t \le s \wedge \strue$.
\end{lem}

\bigskip

For a subset $J \subset [p]$, set
\beq \label{gamma}
\gamma_J = \|\bP_J^\perp \bX \btrue\|. %, \quad \text{and} \quad \gamma_{s,t} = \min_{J \in \cJ_{s,t}} \gamma_J.
\eeq

Assume that both \eqref{inner} and \eqref{quad} hold, which is true with probability at least $1-2 p^{-c}$.  Then, we have that, for all $J\in \mathcal J_{s,t}$:
\begin{eqnarray}
\by^\top (\bP_J - \bP_{\Jtrue}) \by
&\le& - \gamma_J^2 + 2 \gamma_J \sigma \sqrt{(10+2c) (s \vee \strue - t) \log p} + (20+4c) \sigma^2 (s \vee \strue - t) \log p \notag \\
&\le& (40+8c) \sigma^2 (s \vee \strue - t) \log p - \frac12 \gamma_J^2 \label{fidel-ineq-a} \\
%\quad \forall \, J\in \mathcal J_{s,t},\,1\leq s\leq \smax,\, 0\leq t \leq \strue\wedge s.
&\le& (40+8c) \sigma^2 (s \vee \strue - t) \log p. \label{fidel-ineq}
\end{eqnarray}
The first inequality comes from \eqref{fidel}, \eqref{inner} and \eqref{quad}.  The identity $2 a b \le a^2 + b^2$, with $a = \gamma_J/\sqrt{2}$ and $b = \sigma \sqrt{(20+4c) (s \vee \strue - t) \log p}$, justifies the second inequality.

Combining \eqref{ratio} and \eqref{fidel-ineq}, we get
\begin{eqnarray}
\sum_{J\,:\, |J|>[(1+\es )\strue]}^{\smax} \frac{\Pi(J) }{\Pi(\Jtrue)}
&=& \sum_{s =[(1+\es )\strue]}^{\smax} \sum_{t=0}^{\strue} \sum_{J \in \cJ_{s,t}} \frac{{p \choose \strue}}{{p \choose s}} \exp\left(\lambda (\strue - s) + \frac1{2 \sigma^2} \by^\top (\bP_J - \bP_{\Jtrue}) \by \right) \notag\\
&\le& \sum_{s = [(1+\es )\strue]}^{\smax} \sum_{t=0}^{ \strue} \frac{{\strue \choose t} {p - \strue \choose s-t} {p \choose \strue}}{{p \choose s}} \exp\left(\lambda (\strue - s) + (20+4c) (s  - t) \log p  \right), \notag
\end{eqnarray}
where we used the fact that $|\cJ_{s,t}| = {\strue \choose t} {p - \strue \choose s-t}$ in the last inequality.

For the fraction of binomial coefficients, we have
\[
\frac{{\strue \choose t} {p - \strue \choose s-t} {p \choose \strue}}{{p \choose s}} = {s \choose t} {p - s \choose \strue-t}.
\]
%We then use Stirling's bound, which says that for any integer $m$,
%\beq \label{stirling}
%\sqrt{2 \pi} \, m^{m + 1/2} e^{-m} \le m! \le e \, m^{m + 1/2} e^{-m},
%\eeq
%which, for $1 \le k \le m-1$, leads to
%\beq \label{binom}
%\binom{m}{k} \leq \frac{e}{2 \pi} \left(\frac{m}{k}\right)^{k} \left(\frac{m}{m-k}\right)^{m-k + \frac{1}{2}}.
%\eeq
We then use the standard bound on the binomial coefficient
\begin{eqnarray}
\log {s \choose t} + \log {p - s \choose \strue-t}
&\le& (s - t) \log \big(e s/(s-t)\big) + (\strue - t) \log \big(e (p-s)/(\strue - t) \big) \nonumber \\
&\le& 3 (s \vee s^* -t) \log p. \label{inner2}
\end{eqnarray}

Hence, we have so far that
\beq \label{thm-VS-interm1}
\sum_{J\,:\, |J|>[(1+\es )\strue]}^{\smax} \frac{\Pi(J) }{\Pi(\Jtrue)} \le \sum_{s = 0}^{\smax} \sum_{t=0}^{ \strue} \exp\left( A_{s,t} \right),
\eeq
where
\[
A_{s,t} := \omega (s - t) \log p + \lambda (\strue - s), \quad \omega := 23+4c.
\]

Some simple algebra yields
\begin{align}
\sum_{s \ge [(1+\es) \strue]}^{\smax} \sum_{t=0}^{\strue} \exp\left( A_{s,t} \right)
%&=\sum_{s = \strue+1}^{\smax} \sum_{t=0}^{ \strue} \exp\left( A_{s,t} \right)\notag\\
&\le \sum_{s \ge (1+\es) \strue} e^{- (\lambda - \omega \log p) (s -\strue)} \sum_{t=0}^{\strue} e^{(\strue - t) \omega \log p}\notag\\
&\leq \frac{e^{-(\lambda - \omega \log p) \es \strue}}{1-e^{-\lambda + \omega \log p}} \cdot \frac{e^{(\strue+1) \omega \log p}}{e^{\omega \log p} - 1} \label{moment-case} \\
&\leq \frac{p^{-c}}{(1 - p^{-\omega})(1-p^{-c})}, \label{support1}
\end{align}
where we used the fact that $p^\omega \ge 2$, because $p \ge 2$, and also $-(\lambda - \omega \log p) \es \strue + \strue \omega \log p \le - c \log p$, because of \eqref{lambda}.  This shows that
\beqn
\Pi(J: |J| > [(1+\es) \strue]) &\le& \frac{p^{-c}}{(1 - p^{-\omega})(1-p^{-c})} \Pi(\Jtrue) \\
&\le& \frac{p^{-c}}{(1-p^{-c})^2} \Pi(\Jtrue),
\eeqn
using the fact that $\omega \ge c$.  From this, and the fact that $p^{-c} \le 1/2$, we conclude the proof.

%Note that, for any $c > 0$, we have
%\[
%\Pi(J: |J| > [(1+\es) \strue]) \le \frac{p^{-c}}{(1-p^{-c})^2} \Pi(J: |J| \le [(1+\es) \strue]),
%\]
%leading to
%\beq \label{support-alt}
%\Pi(J: |J| \le (1+\es) \strue) \ge (1 - p^{-c})^2 (1 + p^{-c})^{-1} \ge (1 - p^{-c})^3 \ge 1 - 3 p^{-c}.
%\eeq

\subsection{Proof of \thmref{support}}

Let $\nu = \lmin_{(2+\es)\strue}$ for short.  The proof of this result is identical to that of \prpref{support} up to \eqref{fidel-ineq-a}. We now need a lower bound on $\gamma_J$.  For this, we use the following irrepresentability result.
\begin{lem} \label{lem:irrep}
Let $\bX = [\bX_1 \bX_2]$, with smallest singular value $\delta$, and let $\bP_2$ denote the orthogonal projection onto $\bX_2$.  Then for any $\bbeta_1$,
\[
\|(\bI - \bP_2) \bX_1 \bbeta_1\| \ge \delta \|\bbeta_1\|.
\]
\end{lem}

\bigskip
Note that for any $J \in \cJ_{s,t}$ with $s-t \le (1 +\es) \strue$, the smallest singular value of $[\bX_{\Jtrue} \bX_{J \setminus \Jtrue}]$ is bounded from below by $\sqrt{n}\nu$; by \lemref{irrep}, this implies that
\[
\gamma_J = \|(\bI - \bP_{J})(\bX_{\Jtrue}\btrue)\| = \|(\bI - \bP_{J})(\bX_{\Jtrue\setminus J}\bbeta_{\Jtrue\setminus J}^*)\| \geq \sqrt{n}\nu \| \bbeta_{\Jtrue\setminus J}^* \|.
\]
Hence,
\beq \label{gamma1}
\gamma_J \ge  \rho \nu \sqrt{n(\strue - t)}, \quad \forall J\in \mathcal J_{s,t}, \text{ such that } 0\leq  t \leq \strue\wedge s \text{ and } s \le t + (1 +\es) \strue,
\eeq
where we recall that $\rho$ is defined in \eqref{rho}.

In view of \eqref{fidel-ineq-a} and \eqref{gamma1} we have, with probability at least $1-2 p^{-c}$, for all $J\in \mathcal J_{s,t}$
\begin{eqnarray}
\by^\top (\bP_J - \bP_{\Jtrue}) \by
&\le& (40+8c) \sigma^2 (s \vee \strue - t) \log p - \frac12 \gamma_J^2 \notag \\
%\quad \forall \, J\in \mathcal J_{s,t},\,1\leq s\leq \smax,\, 0\leq t \leq \strue\wedge s.
&\le& (40+8c) \sigma^2 (s \vee \strue - t) \log p - \frac12 \rho^2 \nu^2 n (\strue-t) \1_{\{s \le t + (1 +\es) \strue\}}. \label{fidel-ineq-2}
\end{eqnarray}

Next, we have
\begin{eqnarray}
\frac{1}{\Pi(\Jtrue)} &=& \sum_{J\,:\, |J|> [(1+\es)\strue]} \frac{\Pi(J)}{\Pi(\Jtrue)} + \sum_{J\,:\, |J|\le  [(1+\es)\strue]} \frac{\Pi(J)}{\Pi(\Jtrue)}.\label{interm-Pi-1}
\end{eqnarray}
The first sum in the right-hand side was already bounded in \prpref{support}. We concentrate on the second sum.

Combining \eqref{ratio} and \eqref{fidel-ineq-2}, we get
\begin{eqnarray}
\sum_{J\,:\, |J|\le [(1+\es )\strue]} \frac{\Pi(J) }{\Pi(\Jtrue)}
&=& \sum_{s = 0}^{[(1+\es )\strue]} \sum_{t=0}^{s \wedge \strue} \sum_{J \in \cJ_{s,t}} \frac{{p \choose \strue}}{{p \choose s}} \exp\left(\lambda (\strue - s) + \frac1{2 \sigma^2} \by^\top (\bP_J - \bP_{\Jtrue}) \by \right) \notag\\
&\le& \sum_{s = 0}^{[(1+\es )\strue]} \sum_{t=0}^{s \wedge \strue} \frac{{\strue \choose t} {p - \strue \choose s-t} {p \choose \strue}}{{p \choose s}} \exp\big(\lambda (\strue - s) + (20+4c) (s \vee \strue - t) \log p - \eta_{s,t} \big), \notag\\
&\le& \sum_{s = 0}^{[(1+\es )\strue]} \sum_{t=0}^{s \wedge \strue} {s \choose t} {p - s \choose \strue-t} \exp\big(\lambda (\strue - s) + (20+4c) (s \vee \strue - t) \log p - \eta_{s,t} \big), \notag
\end{eqnarray}
where $\eta_{s,t} := \frac1{4\sigma^2} \rho^2 \nu^2 n(\strue-t) \1_{\{s \le t + [(1 +\es) \strue]\}}$.

Next, we use again \eqref{inner2} to get
\beq \label{thm-VS-interm1}
\sum_{J\,:\, |J|\le  [(1+\es)\strue]} \frac{\Pi(J)}{\Pi(\Jtrue)}\le \sum_{s = 0}^{ [(1+\es)\strue]} \sum_{t=0}^{s \wedge \strue} \exp\left( A_{s,t} \right),
\eeq
where
\[
A_{s,t} := \omega (s \vee \strue - t) \log p + \lambda (\strue - s) - \eta_{s,t}, \quad \omega := 23+4c.
\]

Let $\alpha = \frac{\nu^2 n \rho^2}{4\sigma^2} - \omega \log p$, and note that $\alpha \ge 2 \lambda \ge \lambda + c \log p$ by \eqref{lambda} and \eqref{rho}.

When $s \le \strue$, we have $A_{s,t} = -\alpha (s - t) - (\alpha - \lambda) (\strue - s)$, so that
\begin{align}
\sum_{s = 0}^{\strue} \sum_{t=0}^{s} \exp\left( A_{s,t} \right)
%&=\sum_{s = 1}^{\strue-1} \sum_{t=0}^{ s} \exp\left( A_{s,t} \right)\notag\\
&\le \sum_{s = 1}^{\strue} e^{- (\strue - s) c \log p} \sum_{t=0}^{ s} e^{-\alpha (s - t)}\notag\\
%&\leq e^{(-\alpha +\lambda)\strue}\sum_{s = 1}^{\strue}e^{-\lambda s} \frac{e^{\alpha s}-e^{-\alpha }}{1-e^{-\alpha}}\notag\\
&\leq  \frac{1}{(1-e^{-\alpha})(1-p^{-c})}. \label{support2}
\end{align}

%Similarly, when $s = \strue$, we have
%\begin{align*}
%\sum_{t=0}^{\strue} \exp\left( A_{\strue,t} \right)
%%&=\sum_{t=0}^{ \strue} \exp\left( A_{\strue,t} \right)\notag\\
%&= e^{-\alpha \strue} \sum_{t=0}^{ \strue} e^{\alpha t}\notag\\
%&= e^{-\alpha \strue} \frac{e^{\alpha \strue}-e^{-\alpha}}{1-e^{-\alpha}}\notag\\
%&\leq  \frac{1}{1-e^{-\alpha}}\notag\\
%&\leq \frac{e^{-\alpha}}{1-e^{-\alpha}}. %\label{thm-VS-interm3}
%\end{align*}

When $\strue < s \le (1+\es) \strue$, we have $A_{s,t} = -\alpha (\strue - t) - (\lambda - \omega \log p) (s - \strue)$, with $\lambda \ge \omega \log p + c \log p$, leading to
\begin{align}
\sum_{s = \strue+1}^{[(1+\es) \strue]} \sum_{t=0}^{\strue} \exp\left( A_{s,t} \right)
%&=\sum_{s = \strue+1}^{\smax} \sum_{t=0}^{ \strue} \exp\left( A_{s,t} \right)\notag\\
&\le \sum_{s = \strue+1}^{\infty}e^{- (s - \strue) c \log p} \sum_{t=0}^{ \strue} e^{-\alpha (\strue - t)}\notag\\
%&\le e^{[-\frac{\nu\rho^2}{2\sigma^2} +\lambda]\strue}\frac{e^{- (\strue+1) c \log p }}{1-e^{- c \log p}} \frac{e^{\alpha \strue}-e^{-\alpha }}{1-e^{-\alpha}}\notag\\
&\leq \frac{p^{-c}}{(1-e^{-\alpha})(1-p^{-c})}. \label{support3}
\end{align}

Combining \eqref{support1}  with \eqref{interm-Pi-1}-\eqref{support3}, we conclude that
\beqn
\frac1{\Pi(\Jtrue)}
&\leq& \frac{1}{(1-e^{-\alpha})(1-p^{-c})} + \frac{p^{-c}}{(1-e^{-\alpha})(1-p^{-c})} + \frac{p^{-c}}{(1 - p^{-\omega})(1-p^{-c})} \\
&\leq& \frac{1 + 2 p^{-c}}{(1 - p^{-c})^2},
\eeqn
using the fact that $\alpha \ge \omega \ge c$.  From this, we get
\[
\Pi(\Jtrue) \ge (1 - p^{-c})^2 (1 - 2p^{-c}) \ge (1 - 2p^{-c})^2 \ge 1 - 4 p^{-c}.
\]
This concludes the proof of \thmref{support}.  We note that the proof of \eqref{support-refine} is virtually identical.

\subsection{Proof of \thmref{estimation}}

When \eqref{lambda} is satisfied with $\eps \le 1/2$, then $\lambda$ satisfies both the conditions of \prpref{support} and \thmref{prediction}.  Hence, with probability at least $1 - 2 p^{-c} - p^{-c} = 1 - 3 p^{-c}$, we have both that $|\Jmap| \le (1+\eps) \strue$ and \eqref{pred}.  Hence, the support of $\bbmap - \btrue$ is of size at most $(1+\eps) \strue + \strue = (2+\eps) \strue$, and we have
\[
\|\bbmap - \btrue\| \le \frac1{\lmin_{(2+\es)\strue} } \|\bX (\bbmap - \btrue)\|,
\]
with
\[
\|\bX (\bbmap - \btrue)\| = \|\bX \bbmap - \bX \btrue \| \le \sigma \sqrt{8 \strue \lambda},
\]
and the result follows.

\subsection{Proof of \thmref{estimation-1}}

\def\Jsig{J_0}
\def\rhomin{\rho_{\rm min}}
For $r > 0$, we have
\begin{align*}
\P\left(\|\bbmap - \btrue\|_\infty >r \right) &\leq \P\left(\|\bbhat_{\Jtrue} - \btrue\|_\infty >r, \Jmap = \Jtrue \right) +  \P\left(\|\bbmap - \btrue\|_\infty >r, \Jmap \neq \Jtrue \right)\\
&\leq  \P\left(\|\bbhat_{\Jtrue} - \btrue\|_\infty >r\right) +  \P\left(\Jmap \neq \Jtrue \right).
\end{align*}
By \thmref{support}, $\Jmap = \Jtrue$ with probability at least $1 - 2 p^{-c}$, so that the second term on the RHS is bounded by $2 p^{-c}$.
%Let $\Omega$ be the set of $\bz \in \bbR^n$ that satisfy \eqref{inner} and \eqref{quad}.  Letting $\Pi(J; \by)$ emphasize that $\Pi$ depends on $\by$, and remembering that $\E$ denotes the expectation with respect to $\by$, we have
%\begin{eqnarray}
%\P\left(\Jmap \neq \Jtrue \right)
%&=& \E\left[ \Pi(J \neq \Jtrue; \by) \right] \notag\\
%&=& \E\left[ \Pi(J \neq \Jtrue; \by) \1_{\bz \in \Omega} \right] + \E\left[ \Pi(J \neq \Jtrue; \by) \1_{\bz \notin \Omega} \right] \notag\\
%&\le& \E\left[ \Pi(J \neq \Jtrue; \by) \1_{\bz \in \Omega} \right] + \P(\Omega^c) \notag\\
%&\le& 4 p^{-c} + 2 p^{-c} = 6 p^{-c}.\label{sel-var-1}
%\end{eqnarray}
%The first inequality comes from the fact that $\Pi(J \neq \Jtrue; \by) \le 1$.  The second inequality comes from \thmref{support} and its proof, where we learned that $\Pi(J \neq \Jtrue; \by) \le 4 p^{-c}$ on $\Omega$ and $\P(\Omega^c) \le 2 p^{-c}$.

%On the one hand, from \thmref{support} and the union bound,
%\[
%\P\left(\Jhat \neq \Jtrue \right) \le 2 p^{-c} + 4 p^{-c} = 6 p^{-c}.
%\]
%On the other hand

Next, we know that $\bbhat_{\Jtrue} \sim N(\btrue, \sigma^2 \frac{1}{n}\Psi_\star^{-1})$ with $\Psi_\star := \frac{1}{n}\bX_{\Jtrue}^\top \bX_{\Jtrue}$, and in particular, $\widehat{\beta}_{\Jtrue, j} - \beta_{\star, j} \sim \cN(0, \sigma^2 \tau_j^2/n)$, where $\tau_j^2$ is the $j$th diagonal entry of $\Psi_\star^{-1}$.  This matrix being positive semi-definite, its diagonal terms are all bounded from above by its largest eigenvalue, which is the inverse of the smallest eigenvalue of $\Psi_\star$, which in turn is larger than $\nu_{\strue}^{2}$.  Hence, $\Var(\widehat{\beta}_{\Jtrue, j}) \le \sigma^2/(n  \nu_{\strue}^2)$ for all $j \in \Jtrue$, so that a standard tail bound on the normal distribution and the union bound give
\begin{equation}\label{LSE-oracle}
\P\left(\|\bbhat_{\Jtrue} - \btrue\|_\infty >r\right) \le \strue \, \exp\left(-\frac{n \nu_{\strue}^{2} r^2}{2 \sigma^2 }\right).
\end{equation}
Taking $r =  \sigma \sqrt{2 (c+1) \log(p)/(n\nu_{\strue}^2)}$ bounds this by $p^{-c}$, and the desired result follows.

%$r = A\sigma \sqrt{\frac{\log p}{\nu_{\min}(\strue+\smax) n}}$ and

\subsection{Proof of \thmref{estimation-2}}

We have
\begin{eqnarray}
\|\bbmap - \btrue\|_\infty &\le& \sum_{J} \|\bbhat_J - \btrue\|_\infty \Pi(J)\notag\\
&\le& \|\bbhat_{\Jtrue} - \btrue\|_\infty \Pi(\Jtrue) +\sum_{J\neq \Jtrue} \|\bbhat_J - \btrue\|_\infty \Pi(J)\notag\\
&\le&  \|\bbhat_{\Jtrue} - \btrue\|_\infty +\sum_{J\neq \Jtrue} \|\bbhat_J - \btrue\|_\infty \Pi(J).\label{thm:est2-interm-2}
\end{eqnarray}

For any $c>0$, we have with probability at least $1-p^{-c}$, for any $J\subset [p]$ with $\nu_J>0$, that
\begin{eqnarray}
\| \bbhat_{J} \|_\infty &\le& \sqrt{|J|}\| \bbhat_J \| \notag\\
&\le& \frac{\sqrt{|J|}}{\sqrt{n}\nu_J} \|\bX\bbhat_J\|\notag\\
&\le&  \frac{\sqrt{|J|}}{\sqrt{n}\nu_J} \left[ \|\bP_J(\bz)\|  + \|\bP_J^{\perp}(\bX\btrue)\| \right]\notag\\
&\le& \frac{\sqrt{|J|}}{\sqrt{n}\nu_J} \left[\sigma \sqrt{(20+4c) |J|\log p}  + \|\bX\btrue\| \right],\notag
\end{eqnarray}
where we have used Cauchy-Schwarz's inequality in the first line and \eqref{proj} in the last line.

We now assume that $\nu_{\smax}>0$, which implies that $\nu_J>0$ for any $J\subset [p]$ with $|J|\le \smax$.  Combining the previous display with \eqref{LSE-oracle} and \eqref{thm:est2-interm-2} and a union bound argument, we get with probability at least $1-2p^{-c}$, %for any $J\subset [p]$ with $|J|\le \smax$ and $J\neq \Jtrue$, that
\begin{align*}
\|\bbmap - \btrue\|_\infty &\le\sigma\sqrt{\frac{2(c+1)\log p }{ n\nu_{\strue}}}\\ &+ \sum_{J\neq \Jtrue} \left[\frac{\sigma|J|}{\nu_{\smax}}\sqrt{(20+4c)\log p}  + \frac{\sqrt{|J|}}{\sqrt{n}\nu_{\smax}}\|\bX\btrue\| + \|\btrue\|_\infty\right]\Pi(J).
\end{align*}
Next, we combine the above display with \eqref{support-refine} and a union bound argument to get with probability at least $1-4p^{-c}$ that
\begin{eqnarray*}
\|\bbmap- \btrue\|_\infty &\le& \sigma\sqrt{\frac{2(c+1)\log p }{ n\nu_{\strue}}} + \frac{4 p^{-c}}{\nu_{\smax}} \left[  \sigma\sqrt{(20+4c)\frac{\log p}{n}}  + \frac{\|\bX\btrue\|}{\sqrt{n}} + \nu_{\smax}\|\btrue\|_\infty \right].
\end{eqnarray*}
Note that the same reasoning applied to $\widetilde \bbeta$ yields the same $l_\infty$-norm estimation bound with $\nu_{\smax}$ replaced by $\nu_{\min}$.

We now assume that $\nu_{\strue+\smax} > 0$. Then, for any $J\subset [p]$ with $|J|\le \smax$, we have
$$
\|\bbhat_J - \btrue\|_\infty \le  \| \bbhat_J - \btrue\| \le \frac{\|\bX \bbhat_J - \bX\btrue\|}{\sqrt{n}\nu_{\strue+\smax}}.
$$
%Recall that $\cJ = \{J: \|\bxi_J\| > \sigma \sqrt{10 \strue \lambda}\}$
% and set $\tilde\cJ = \left\{ J\subset [p]\,:\, |J|\le \smax,\, J\neq \Jtrue\right\}$.
 Combining this last inequality with \eqref{thm:est2-interm-2}, we get
\begin{eqnarray*}
\|\bbmap - \btrue\|_\infty &\le&  \|\bbhat_{\Jtrue} - \btrue\|_\infty + \frac{1}{\sqrt{n}\nu_{\strue+\smax}} \sum_{J\notin \cJ, J \neq \Jtrue} \|\bxi_J\|\Pi(J) + \frac{1}{\sqrt{n}\nu_{\strue+\smax}} \sum_{J\in \cJ} \|\bxi_J\|\Pi(J)\\
&\le& \|\bbhat_{\Jtrue} - \btrue\|_\infty + \frac{\sigma \sqrt{10 \strue}}{\sqrt{n}\nu_{\strue+\smax}} \Pi(\cJ^{c} \setminus \Jtrue) + \frac{1}{\sqrt{n}\nu_{\strue+\smax}} \sum_{J\in \cJ} \|\bxi_J\|\Pi(J),
\end{eqnarray*}
where we recall that $\bxi_J= \bX\bbhat_{J} - \bX\btrue$ and $\cJ = \left\{J\subset [p]\,:\, \|\bxi_J\|>\sigma \sqrt{10\strue \lambda}\right\}$.
In view of \thmref{support}, we have with probability at least $1-2p^{-c}$ that
$$
\Pi(\cJ^{c} \setminus \Jtrue) \le 1 - \Pi(\Jtrue) \le 4 p^{-c};
$$
and in view of \eqref{sum-cJ},
$$
\sum_{J\in \cJ} \|\bxi_J\|\Pi(J)
\le 2\sqrt{10}\sigma p^{-\strue}.
$$
Combining the three last displays with \eqref{LSE-oracle}, we get the result.

\subsection{Proofs of auxiliary results}

\lemref{proj} is a special case of \lemref{quad} where $\Jtrue = \emptyset$, and we prove \lemref{quad} below.

\subsubsection{Proof of \lemref{inner}}
First, note that $u_J := \< \bP_J^\perp \bX \btrue, \bz \> \sim \cN(0, \sigma^2 \gamma_J^2)$, where $\gamma_J$ is defined in \eqref{gamma}, so that $v_J := u_J/(\sigma \gamma_J) \sim \cN(0,1)$.  By the union bound and a standard tail bound on the normal distribution, for $a > 0$, we have
\beqn
\pr{\max_{J \in \cJ_{s,t}} v_J^2 > a^2}
&\le& {\strue \choose t} {p - \strue \choose s-t} \exp(-a^2/2).
\eeqn
As in \eqref{inner2}, we have
\begin{eqnarray}
\log {\strue \choose t} + \log {p - \strue \choose s-t}
&\le& (\strue - t) \log (e \strue) + (s - t) \log (e p) \nonumber \\
&\le& 3 (s \vee \strue -t) \log p. \label{inner3}
\end{eqnarray}
Hence,
\[
\pr{\max_{J \in \cJ_{s,t}} v_J^2 > (10+2c) (s \vee \strue -t) \log p} \le \exp\big(-(2+c) (s \vee \strue -t) \log p \big) \le p^{-(2+c)},
\]
since $s \vee \strue - t = 0$ would imply $J = \Jtrue$.  We then apply the union bound again,
\[
\pr{\max_{s,t} \max_{J \in \cJ_{s,t}} \frac{v_J^2}{s \vee \strue -t}  > (10+2c) \sigma^2 \log p} \le \smax \, (s \wedge \strue + 1) p^{-(2+c)} \le p^{-c},
\]
which the result we wanted.

\subsubsection{Proof of \lemref{quad}}
Fix $J \in \cJ_{s,t}$.
First, we notice that
\[
\bz^\top (\bP_J - \bP_{\Jtrue}) \bz = \bz^\top (\bP_J - \bP_{J \cap \Jtrue}) \bz - \bz^\top (\bP_{\Jtrue} - \bP_{J \cap \Jtrue}) \bz \le \bz^\top (\bP_J - \bP_{J \cap \Jtrue}) \bz,
\]
since $\bP_{\Jtrue} - \bP_{J \cap \Jtrue}$ is an orthogonal projection, and therefore positive semidefinite.  And $\bQ_J := \bP_J - \bP_{J \cap \Jtrue}$ is also an orthogonal projection, of rank $s - t$, so that $\|\bQ_J \bz\|^2 \sim \sigma^2 \chi_{s-t}^2$.
Chernoff's Bound applied to the chi-square distribution yields
\[\log \pr{\chi^2_m > a} \leq - \frac{m}{2} (a/m - 1 - \log (a/m)) \leq -\frac{a}4, \quad \forall a \geq 2 m.\]
The union bound and \eqref{inner3}, and this tail bound, yields
\[
\pr{\max_{J \in \cJ_{s,t}} \|\bQ_J \bz\|^2 > (20+4c) \sigma^2 (s \vee \strue-t) \log p} \le \exp\left(-(2+c)  (s \vee \strue -t) \log p\right).
\]
The rest of the proof is exactly the same as that of \lemref{inner}.

\subsection{An irrepresentability result}

We have
\beqn
\|(\bI - \bP_2) \bX_1 \bbeta_1\|^2
&=& \min_{\bbeta_2} \|\bX_1 \bbeta_1 + \bX_2 \bbeta_2\|^2 \\
&=& \min_{\bbeta_2}  \bbeta \bX^\top \bX \bbeta \\
&\ge& \min_{\bbeta_2}  \delta^2 \|\bbeta\|^2 \\
&=& \delta^2 \|\bbeta_1\|^2,
\eeqn
where $\bbeta := (\bbeta_1, \bbeta_2)$, implying $\|\bbeta\|^2 = \|\bbeta_1\|^2 + \|\bbeta_2\|^2$.

\subsection*{Acknowledgements}
We would like to thank Pierre Alquier, Philippe Rigollet and Alexander Tsybakov for sharing their R codes with us.
This work was partially supported by XXX.

% bibliograhy
\renewcommand{\refname}{\begin{center}{\normalfont \huge References}\end{center}}
\bibliographystyle{chicago}
\bibliography{biblio-karim}

\end{document}